\message{by Udo Hertrich-Jeromin, August/September 1997}\immediate\write1{}

\magnification = \magstep0
\vsize = 525dd
\hsize = 27cc
\topskip = 13dd
\hoffset = 1.8cm
\voffset = 1.8cm

\message{more fonts: petit and fraktur;}

\font\bfbig = cmbx10 scaled \magstep2   

\font\eightrm = cmr10 scaled 800        
\font\sixrm = cmr7 scaled 850
\font\fiverm = cmr5
\font\eighti = cmmi10 scaled 800
\font\sixi = cmmi7 scaled 850
\font\fivei = cmmi5
\font\eightit = cmti10 scaled 800
\font\eightsy = cmsy10 scaled 800
\font\sixsy = cmsy7 scaled 850
\font\fivesy = cmsy5
\font\eightsl = cmsl10 scaled 800
\font\eighttt = cmtt10 scaled 800
\font\eightbf = cmbx10 scaled 800
\font\sixbf = cmbx7 scaled 850
\font\fivebf = cmbx5

\font\fivefk = eufm5                    
\font\sixfk = eufm7 scaled 850
\font\sevenfk = eufm7
\font\eightfk = eufm10 scaled 800
\font\tenfk = eufm10

\newfam\fkfam
        \textfont\fkfam=\tenfk \scriptfont\fkfam=\sevenfk
                \scriptscriptfont\fkfam=\fivefk

\def\eightpoint{%
        \textfont0=\eightrm \scriptfont0=\sixrm \scriptscriptfont0=\fiverm
                \def\rm{\fam0\eightrm}%
        \textfont1=\eighti  \scriptfont1=\sixi  \scriptscriptfont1=\fivei
                \def\oldstyle{\fam1\eighti}%
        \textfont2=\eightsy \scriptfont2=\sixsy \scriptscriptfont2=\fivesy
        \textfont\itfam=\eightit \def\it{\fam\itfam\eightit}%
        \textfont\slfam=\eightsl \def\sl{\fam\slfam\eightsl}%
        \textfont\ttfam=\eighttt \def\tt{\fam\ttfam\eighttt}%
        \textfont\bffam=\eightbf \scriptfont\bffam=\sixbf
                \scriptscriptfont\bffam=\fivebf \def\bf{\fam\bffam\eightbf}%
        \textfont\fkfam=\eightfk \scriptfont\fkfam=\sixfk
                \scriptscriptfont\fkfam=\fivefk
        \rm}
\skewchar\eighti='177\skewchar\sixi='177\skewchar\eightsy='60\skewchar\sixsy='60

\def\fk{\fam\fkfam}
\def\petit{\eightpoint}%


\def\today{\ifcase\month\or
        January\or February\or March\or April\or May\or June\or
        July\or August\or September\or October\or November\or December\fi
        \space\number\day, \number\year}
\def\newline{\hfil\break}

\message{ lay out: title page,}

\newdimen\fullhsize \fullhsize=\hsize

\newbox\TITLEBOX \setbox\TITLEBOX=\vbox{}
\newdimen\TITLEHEIGHT
\long\def\Title#1#2#3#4{\setbox\TITLEBOX=\vbox{%
        \topglue3truecm%
        \noindent{\bfbig#1}\vskip12pt%
        \noindent{\bf#2}\vskip6pt%
        \noindent{\petit#3}\vskip10pt%
        \noindent{\petit\today}\vskip32pt%
        \noindent{\bf Summary.}\enspace#4\vskip32pt
        \relax}%
        \TITLEHEIGHT=\ht\TITLEBOX%
        \if\IFTHANKS Y\advance\TITLEHEIGHT by\ht\THANKBOX\fi%
        \advance\vsize by-\TITLEHEIGHT%
        \relax}

\newtoks\RunAuthor\RunAuthor={} \newtoks\RunTitle\RunTitle={}
\def\ShortTitle#1#2{\RunTitle={#1}\RunAuthor={#2}}
\headline={\ifnum\pageno=1 {\hfil} \else
        \ifodd\pageno {\petit{\the\RunTitle}\hfil\folio}
        \else {\petit\folio\hfil{\the\RunAuthor}} \fi \fi}
\def\makeheadline{\vbox{%
        \hbox to\fullhsize{\the\headline}%
        \vss\nointerlineskip\kern2pt%
        \hbox to\fullhsize{\hrulefill}\kern7pt}}

\newbox\THANKBOX \let\IFTHANKS N
\setbox\THANKBOX=\vbox{\kern12pt\hbox to\fullhsize{\hrulefill}\kern4pt}
\def\Thanks#1#2{\nobreak${}^{#1}$\global\let\IFTHANKS Y%
        \global\setbox\THANKBOX=\vbox{\parindent20pt\baselineskip9pt%
        \unvbox\THANKBOX{\petit\item{${}^{#1}$}{#2}} }}
\footline={\ifnum\pageno=1 \if\IFTHANKS Y\box\THANKBOX\fi
        \else\hfill\fi}
\def\makefootline{\hbox to\fullhsize{\the\footline}}


\def\OneColumn{\output={%
        \shipout\vbox{\makeheadline%
                \ifnum\pageno=1\hbox to\fullhsize{\box\TITLEBOX\hfil}\fi%
        \hbox to\fullhsize{\pagebody}%
        \makefootline}%
                \ifnum\pageno=1\global\advance\vsize by\TITLEHEIGHT\fi%
        \advancepageno%
        \ifnum\outputpenalty>-20000%
        \else\dosupereject%
        \fi}%
        }%
\OneColumn

\newbox\LEFTCOLUMN
\def\TwoColumn{\let\lr=L \hsize=.485\fullhsize \output={%
        \if L\lr\global\setbox\LEFTCOLUMN=\leftline{\pagebody}\global\let\lr=R%
        \else \shipout\vbox{\makeheadline%
                \ifnum\pageno=1\hbox to\fullhsize{\box\TITLEBOX\hfil}\fi%
        \hbox to\fullhsize{\box\LEFTCOLUMN\hfil\leftline{\pagebody}}%
        \makefootline}%
                \ifnum\pageno=1\global\advance\vsize by\TITLEHEIGHT\fi%
        \advancepageno\global\let\lr=L%
        \fi%
        \ifnum\outputpenalty>-20000%
        \else\dosupereject%
        \fi}%
        }%


\def\PLabel#1{\xdef#1{\nobreak(p.\the\pageno)}}

\message{headings and}

\newcount\SECNO \SECNO=0
\newcount\SUBSECNO \SUBSECNO=0
\newcount\SUBSUBSECNO \SUBSUBSECNO=0
\def\Section#1{\SUBSECNO=0\SUBSUBSECNO=0 \advance\SECNO by 1
        \goodbreak\vskip14pt\noindent{\bf\the\SECNO .\ #1}
        \gdef\Label##1{\xdef##1{\nobreak\the\SECNO}}
        \nobreak\vskip8pt\noindent\kern0pt}
\def\SubSection#1{\SUBSUBSECNO=0 \advance\SUBSECNO by 1
        \goodbreak\vskip14pt\noindent{\it\the\SECNO.\the\SUBSECNO\ #1}
        \gdef\Label##1{\xdef##1{\nobreak\the\SECNO.\the\SUBSECNO}}
        \nobreak\vskip8pt\noindent\kern0pt}
\def\SubSubSection#1{\advance\SUBSUBSECNO by 1
        \goodbreak\vskip14pt\noindent{\rm\the\SECNO.\the\SUBSECNO.\the\SUBSUBSECNO\ #1}
        \gdef\Label##1{\xdef##1{\nobreak\the\SECNO.\the\SUBSECNO.\the\SUBSUBSECNO}}
        \nobreak\vskip8pt\noindent\kern0pt}

\message{definitions;}

\long\def\Definition#1#2{\medbreak\noindent{\bf Definition%
        #1.\enspace}{\it#2}\medbreak\smallskip\relax}
\long\def\Theorem#1#2{\medbreak\noindent{\bf Theorem%
        #1.\enspace}{\it#2}\medbreak\smallskip\relax}
\long\def\Lemma#1#2{\medbreak\noindent{\bf Lemma%
        #1.\enspace}{\it#2}\medbreak\smallskip\relax}


\newcount\FOOTNO \FOOTNO=0
\long\def\Footnote#1{\global\advance\FOOTNO by 1
        {\parindent=20pt\baselineskip=9pt%
        \footnote{\nobreak${}^{\the\FOOTNO)}$}{\petit#1\par\vskip-9pt}%
        }\gdef\Label##1{\xdef##1{\nobreak\the\FOOTNO}}}

\message{referencing system: equation numbers,}

\newcount\EQNO \EQNO=0
\def\Eqno{\global\advance\EQNO by 1 \eqno(\the\EQNO)%
        \gdef\Label##1{\xdef##1{\nobreak(\the\EQNO)}}}

\message{figures and tables and}

\newcount\FIGNO \FIGNO=0
\def\Fcaption#1{\global\advance\FIGNO by 1
        {\petit{\bf Fig. \the\FIGNO.~}#1}
        \gdef\Label##1{\xdef##1{\nobreak\the\FIGNO}}}

\newcount\TABNO \TABNO=0
\def\Tcaption#1{\global\advance\TABNO by 1
   {\petit{\bf Table. \the\TABNO.~}#1}
   \gdef\Label##1{\xdef##1{\nobreak\the\TABNO}}}

\message{the bibliography;}

\newcount\REFNO \REFNO=0
\newbox\REFBOX \setbox\REFBOX=\vbox{}
\def\BegRefs{\message{reading the references}\setbox\REFBOX\vbox\bgroup
        \parindent18pt\baselineskip9pt\petit}
\def\Ref#1{\message{.}\global\advance\REFNO by 1 \ifnum\REFNO>1\vskip3pt\fi
        \item{\the\REFNO .~}\xdef#1{\nobreak[\the\REFNO]}}
\def\EndRefs{\par\egroup\message{done}}
\def\References{\goodbreak\vskip21pt\leftline{\bf References}
        \nobreak\vskip12pt\unvbox\REFBOX\vskip21pt\relax}

\message{finally: some math things.  }

\def\N{I\kern-.8ex N} \def\Z{\raise.72ex\hbox{${}_{\not}$}\kern-.45ex {\rm Z}}
\def\Q{\raise.82ex\hbox{${}_/$}\kern-1.35ex Q} \def\R{I\kern-.8ex R}
\def\C{\raise.87ex\hbox{${}_/$}\kern-1.35ex C} \def\H{I\kern-.8ex H}
\def\HP{\H\!\!P}

\def\Im{{\rm Im}} \def\Re{{\rm Re}}
\def\det{{\rm det}} \def\tr{{\rm tr}} \def\SD{{\cal D}} \def\P{{\cal P}}



\Title{Supplement on Curved flats in the space of point pairs and
        Isothermic surfaces: A Quaternionic Calculus}
        {Udo Hertrich-Jeromin\Thanks{\ast}{Partially supported by the
        Alexander von Humboldt Stiftung and NSF Grant DMS93-12087.}}
        {Dept.~Math.~\&~Stat., GANG, University of Massachusetts,
        Amherst, MA 01003}
{A quaternionic calculus for surface pairs in the conformal 4-sphere is
elaborated. This calculus is then used to discuss the relation between
curved flats in the symmetric space of point pairs and Darboux and
Christoffel pairs of isothermic surfaces. A new viewpoint on relations
between surfaces of constant mean curvature in certain space forms is
presented --- in particular, a new form of Bryant's Weierstrass type
representation for surfaces of constant mean curvature 1 in hyperbolic
3-space is given.}

\ShortTitle{U.Jeromin}{Curved flats and Isothermic surfaces}

\BegRefs
\Ref\Asla H.~Aslaksen: {\it Quaternionic Determinants\/};
        Math.\ Intell.\ {\bf 18.3} (1996) 57-65
\Ref\Blaschke W.~Blaschke: {\it Vorlesungen \"uber Differentialgeometrie
        III\/}; Springer, Berlin 1929
\Ref\Bryant R.~Bryant: {\it Surfaces of mean curvature one in hyperbolic
        space\/}; Ast\'{e}risque {\bf 154-155} (1987) 321-347
\Ref\BJPP F.~Burstall, U.~Hertrich-Jeromin, F.~Pedit, U.~Pinkall:
        {\it Isothermic surfaces and Curved flats\/}; Math.~Z. {\bf 225}
        (1997) 199-209
\Ref\FePe D.~Ferus, F.~Pedit: {\it Curved flats in Symmetric spaces\/};
        Manuscripta Math.\ {\bf 91} (1996) 445-454
\Ref\HHP U.~Hertrich-Jeromin, T.~Hoffmann, U.~Pinkall: {\it A discrete
        version of the Darboux transform for isothermic surfaces\/};
        to appear in A.~Bobenko, R.~Seiler, {\it Discrete integrable
        Geometry and Physics\/}, Oxford Univ.\ Press, Oxford 1997
\Ref\JePe U.~Hertrich-Jeromin, F.~Pedit: {\it Remarks on the Darboux
        transform of isothermic surfaces\/}; to appear in Documenta Math.
\Ref\Palmer B.~Palmer: {\it Isothermic surfaces and the Gauss map\/};
        Proc.\ Amer.\ Math.\ Soc.\ {\bf 104} (1988) 876-884
\Ref\Study E.~Study: {\it Ein Seitenst\"uck zur Theorie der linearen
        Transformationen einer komplexen Ver\-\"anderlichen, Teile I-IV\/};
        Math.~Z. {\bf 18} (1923) 55-86, 201-229 and {\bf 21} (1924)
        45-71, 174-194
\Ref\Wilker J.~Wilker: {\it The Quaternion formalism for M\"obius groups
        in four or fewer dimensions\/}; Lin.\ Alg.\ Appl.\ {\bf 190} (1993)
        99-136
\Ref\UmYa M.~Umehara, K.~Yamada: {\it A parametrization of the Weierstrass
        formulae and perturbation of complete minimal surfaces in $\R^3$
        into the hyperbolic 3-space\/}; J.~reine angew.\ Math.\ {\bf 432}
        (1992) 93-116
\EndRefs

\Section{Introduction}
It is well known that the orientation preserving M\"obius transformations of
the ``conformal 2-sphere'' $S^2\cong\C\cup\{\infty\}$ can be described as
fractional linear transformations $z\mapsto{a_{11}z+a_{12}\over a_{21}z+a_{22}}$
where $a=(a_{ij})\in Sl(2,\C)$.
The reason for this fact is that the conformal 2-sphere $S^2\cong\C P^1$
can be identified with the complex projective line.
Introducing homogeneous coordinates $p=v_p\C$, $v_p\in\C^2$, on $\C P^1$
the special linear group $Sl(2,\C)$ acts on $\C P^1$ by M\"obius transformations
--- which are, for 1-dimensional projective spaces, identical with
projective transformations --- via $v_p\C\mapsto Av_p\C=v_q\C$.
Thus, in affine coordinates one has
$$\left(\matrix{z\cr1\cr}\right)\mapsto
        \left(\matrix{a_{11}&a_{12}\cr a_{21}&a_{22}\cr}\right)\cdot
        \left(\matrix{z\cr1\cr}\right)\simeq
        \left(\matrix{{a_{11}z+a_{12}\over a_{21}z+a_{22}}\cr1\cr}\right)$$
This (algebraic) model of M\"obius geometry in dimension 2 complements the
(``geometric'') model commonly used in differential geometry: here, the
conformal $2$-sphere (or, more general, the conformal $n$-sphere) is
considered as a quadric in the real projective 3-space $\R P^3$ and
the group of M\"obius transformations is isomorphic to the group of
projective transformations of $\R P^3$ that map the ``absolute quadric''
$S^2$ onto itself (cf.\Blaschke).
Equipping the space of homogeneous coordinates of $\R P^3$ with a Lorentz
scalar product that has the points of $S^2$ as isotropic (null) lines,
the M\"obius group can be identified with the pseudo orthogonal group of
this Minkowski space $\R^4_1$.

Several attempts have been made to generalize the described algebraic model to
higher dimensions --- in particular to dimensions 3 and 4, by using quaternions
(cf.\Study,\Wilker): analogous to the above model, the conformal 4-sphere is
identified with the quaternionic projective line, $S^4\cong\HP^1$, with
$Sl(2,\H)$ acting on it by M\"obius transformations.
In order to use such an ``algebraic model'' in M\"obius {\it differential}
geometry, it is not enough to describe the underlying space and the M\"obius
group acting on it. One also needs a convenient description for (hyper-)
spheres since the geometry of surfaces in M\"obius geometry is often closely
related to the geometry of an enveloped sphere congruence (cf.\Blaschke).
For example, Willmore surfaces in $S^3$ can be related to minimal surfaces in
the space of 2-spheres in $S^3$, and the geometry of isothermic surfaces is
related to that of ``sphere surfaces'' with flat normal bundle, ``Ribaucour
sphere congruences''.

One way is to identify a hypersphere $s\subset\HP^1$ with the inversion at
this sphere. The problem with this approach is, that only the orientation
preserving M\"obius transformations are naturally described in the algebraic
model --- but, inversions are orientation reversing M\"obius transformations.
Adjoining the (quaternionic) conjugation as a basic orientation reversing
M\"obius transformation and working with the larger group of {\it all} M\"obius
transformations, works relatively fine for 2-dimensional M\"obius geometry,
but turns into a nightmare\Footnote{Identifying 2-spheres in $S^3\subset
S^4\cong\HP^1$ with inversions in $S^4$ provides a solution: as the composition
of two inversions at hyperspheres, the inversion at a 2-sphere in $S^4$ is
orientation preserving.} in dimension 4 since the quaternions form a non
commutative field.

Another way is to identify a sphere $s\subset S^4\cong\HP^1$ with that
quaternionic hermitian form on the space $\H^2$ of homogeneous coordinates
that has this sphere $s$ as a null cone.
After discussing some basics in quaternionic linear algebra we will follow
this approach --- to obtain not only a description for the space of spheres
but also to establish the relation with the classical ``geometric'' model of
M\"obius geometry: the space of quaternionic hermitian forms will canonically
turn into a real six dimensional Minkowski space, the classical model space.

Using this setup, we discuss the geometry of surface pairs, maps into the
symmetric space of point pairs in $\HP^1$.  In M\"obius differential geometry,
surface {\it pairs} occur in various situations: in the context of Willmore
surfaces (the dual) as well as in the context of isothermic surfaces
(Christoffel and Darboux pairs).
The latter will be examined in the remaining part of the paper, on one side
to see the calculus at work, on the other side to demonstrate some new results:
although the relation between Darboux pairs of isothermic surfaces in $S^3$ and
curved flats in the space of point pairs was already established in \BJPP\
it might be of interest to see that this relation also holds in the higher
codimension case of Darboux pairs in $\HP^1$ (cf.\Palmer).
Also, our quaternionic calculus provides very elegant characterizations for
Darboux and Christoffel pairs of isothermic surfaces that led to the discovery
of the Riccati type equation (cf.\JePe) for the Darboux transformation of
isothermic surfaces, and hence was crucial for the definition of the discrete
version of the Darboux transformation and the (geometric) definition of
discrete cmc nets (cf.\HHP).

In the last section, we study minimal and constant mean curvature surfaces
in 3-dimensional spaces of constant curvature. These are ``special''
isothermic surfaces, and a suitable Christoffel transform in $\R^3$ can
be determined algebraically (in the general case, an integration has to be
carried out).
Examining the effect of the spectral parameter that comes with a curved flat,
we obtain a new interpretation for the relations between surfaces of constant
curvature in certain space forms.
For example, the well known relation between minimal surfaces in the (metric)
3-sphere and surfaces of constant mean curvature in Euclidean space, as well
as the relation between minimal surfaces in Euclidean 3-space and surfaces of
constant mean curvature 1 in hyperbolic 3-space are discussed.
In case of the constant mean curvature 1 surfaces in hyperbolic 3-space, a new
form of Bryant's Weierstrass type representation \Bryant\ is given.
In this context, the classical Weierstrass representation for minimal surfaces
in Euclidean 3-space is described as a Goursat type transform of the plane ---
similar to the way certain surfaces of constant Gauss curvature are described
as a B\"acklund transform of a line.
In fact, the classical Goursat transformation for minimal surfaces is
generalized for isothermic surfaces in Euclidean space.

\Section{The Study determinant}
Throughout this paper we will use various well known models \Asla\ for the
non commutative field of quaternions:
$$\matrix{\H&\cong&\{a+v\,|\,a\in\R\cong\Re\H,v\in\R^3\cong\Im\H\}\hfill\cr
        &\cong&\{a_0+a_1i+a_2j+a_3k\,|\,a_0,a_1,a_2,a_3\in\R\}\hfill\cr
        &\cong&\{x+y\,j\,|\,x,y\in\C\}\hfill\cr
        &\cong&\{A\in M(2\times2,\C)\,|\,\tr A\in\R,A+A^{\ast}\in\R I\}.\cr}$$
Herein, we can identify $i,j,k$ with the standard basis vectors of $\R^3
\cong\Im\H$: if $v,w\in\Im\H$ are two ``vectors'' their product $v\,w=
-v\cdot w+v\times w$ which is equivalent to the identities $i^2=j^2=k^2=-1$,
$ij=k=-ji$, $jk=i=-kj$ and $ki=j=-ik$. Obviously, the first model will turn
out particularly useful when focussing on the geometry of 3-space while the
decomposition $\H\cong\C+\C\,j$ will prove useful in the context of surfaces,
2-dimensional submanifolds, since their tangent planes (and normal planes)
carry a natural complex structure. We will switch between these models as it
appears convenient.

As the quaternions can be thought of as a Euclidean 4-space, $\R^4\cong\H$,
the (conformal) 4-sphere $S^4\cong\R^4\cup\{\infty\}$ can be identified with
the quaternionic projective line: $S^4\cong\HP^1=\{\hbox{lines through 0 in
$\H^2$}\}$. Thus, a point $p\in S^4$ of the conformal 4-sphere is described
by its homogeneous coordinates $v_p\in\H^2$; and its stereographic projection
onto Euclidean 4-space $\R^4\cong\{v\in\H^2\,|\,v_2=1\}$ is obtained by
normalizing the second component of $v_p$.

Since the quaternions form a non commutative field, we have to agree whether
the scalar multiplication in a quaternionic vector space is from the right or
left: in this paper, $\H^2$ will be considered a {\it right vector space} over
the quaternions. This way, quaternionic linear transformations can be
described by the multiplication (of column vectors) with (quaternionic)
matrices from the {\it left\/}: $A(v\lambda)=(Av)\lambda$. For a quaternionic
2-by-2 matrix $A\in M(2\times2,\H)$ we introduce the Study determinant%
\Footnote{Note, that the notion of determinant makes sense for self adjoint
matrices $A\in M(2\times2,\H)$.} \Asla\ (cf.~Study's ``Nablafunktion'' \Study)
$$\matrix{\SD(A)&:=&\det(A^{\ast}A)\hfill\cr
        &\hfill=&|a_{11}|^2|a_{22}|^2+|a_{12}|^2|a_{21}|^2
          -(\bar{a}_{11}a_{12}\bar{a}_{22}a_{21}+
            \bar{a}_{21}a_{22}\bar{a}_{12}a_{11}).\cr}$$
This is exactly the determinant of the complex 4-by-4 matrix corresponding
to $A$ when using the complex matrix model for the quaternions. Thus, $\SD$
clearly satisfies the usual multiplication law, $\SD(AB)=\SD(A)\SD(B)$, and
vanishes exactly when $A$ is singular. The multiplication law implies that
$\SD$ is actually an invariant of the linear transformation described by a
matrix: $\SD(U^{-1}A\,U)=\SD(A)$ for any basis transformation $U:\H^2\to\H^2$.
Also note that $0\leq\SD(A)\in\R$.

\Definition{}{The general and special linear groups of $\H^2$ will be denoted
by $$\matrix{Gl(2,\H)&:=&\{A\in M(2\times2,\H)\,|\,\SD(A)\neq0\}\hfill\cr
        Sl(2,\H)&:=&\{A\in M(2\times2,\H)\,|\,\SD(A)=1\}.\cr}$$}

\noindent With the help of Study's determinant, the inverse of a quaternionic
2-by-2 matrix $A\in Gl(2,\H)$ can be expressed directly as
$$A^{-1}={1\over\SD(A)}\left(\matrix{
        |a_{22}|^2\bar{a}_{11}-\bar{a}_{21}a_{22}\bar{a}_{12}&
        |a_{12}|^2\bar{a}_{21}-\bar{a}_{11}a_{12}\bar{a}_{22}\cr
        |a_{21}|^2\bar{a}_{12}-\bar{a}_{22}a_{21}\bar{a}_{11}&
        |a_{11}|^2\bar{a}_{22}-\bar{a}_{12}a_{11}\bar{a}_{21}\cr}\right).$$
Note also, that $Sl(2,\H)$ is a 15-dimensional Lie group --- it will turn out
to be a double cover of the identity component of the M\"obius group of $S^4$.

Considering $\SD:\H^2\times\H^2\to\R$ as a function of two (column) vectors
we see that $\SD(v,v+w)=\SD(v,w)$ and $\SD(v,w\lambda)=|\lambda|^2\SD(v,w)$
--- similar formulas holding for the first entry since $\SD$ is symmetric:
$\SD(v,w)=\SD(w,v)$. Reformulating our previous statement, we also obtain
that $\SD(v,w)=0$ if and only if $v$ and $w$ are linearly dependent%
\Footnote{All these properties are also easily checked directly, without
using the complex matrix representation of the quaternions.}.
Particularly, if $v$ and $w$ are points in an affine quaternionic line,
say the Euclidean 4-space $\{v\in\H^2\,|\,v_2=1\}$, then
$\SD(v,w)=|v_1-w_1|^2$ measures the distance between $v$ and $w$ with
respect to a Euclidean metric. This fact can be used to express the
cross ratio of four points in Euclidean 4-space (cf.\HHP) in terms of the
Study determinant\Footnote{For a more complete discussion of the {\it complex}
cross ratio of four points in space consult \HHP.}:
$$|D\!V(h_1,h_2,h_3,h_4)|^2={
        \SD{\petit\left(\matrix{h_1&h_2\cr1&1\cr}\right)}
        \SD{\petit\left(\matrix{h_3&h_4\cr1&1\cr}\right)}\over
        \SD{\petit\left(\matrix{h_2&h_3\cr1&1\cr}\right)}
        \SD{\petit\left(\matrix{h_4&h_1\cr1&1\cr}\right)}}.$$
The expression on the right hand is obviously invariant under individual
rescalings of the vectors which shows that the ross ratio is, in fact, an
invariant of four points in the quaternionic projective line $\HP^1$.

\Section{Quaternionic hermitian forms}
will be a key tool in our calculus for M\"obius geometry:
any quaternionic hermitian form $s:\H^2\times\H^2\to\H$,
$$\matrix{s(v,w_1\lambda+w_2\mu)&=&s(v,w_1)\lambda+s(v,w_2)\mu\cr
        s(v_1\lambda+v_2\mu,w)&=&\bar{\lambda}s(v_1,w)+\bar{\mu}s(v_2,w)\cr
        \hfill          s(w,v)&=&\overline{s(v,w)},     \hfill\cr}$$
is determined by its values on a basis $(e_1,e_2)$ of $\H^2$,
$s_{ij}=s(e_i,e_j)$. Since $s$ is hermitian, $s_{11},s_{22}\in\R$
and $s_{21}=\bar{s}_{12}\in\H$, the quaternionic hermitian forms on
$\H^2$ form a 6-dimensional (real) vector space.
Clearly, $Gl(2,\H)$ operates on this vector space via $(A,s)\mapsto As
:=[(v,w)\mapsto s(Av,Aw)]$, or, in the matrix representation of $s$, via
$(A,s)\mapsto A^{\ast}sA$. A straightforward calculation shows that
$\det(As)=\SD(A)\det(s)$. This enables us to introduce a Lorentz scalar product
$$\langle s,s\rangle:=-\det(s)=|s_{12}|^2-s_{11}s_{22}$$
on the space $\R^6_1$ of quaternionic hermitian forms, which is well defined
up to a scale\Footnote{At this point, we notice that the geometrically
significant space is the projective 5-space $\R\!P^5$ with absolute quadric
$Q=\{\R x\,|\,\langle x,x\rangle=0\}$, not its space of homogeneous coordinates,
$\R^6_1$.} (or, the choice of a basis in $\H^2$). Fixing a scaling of this
Lorentz product, the special linear transformations act as isometries on
$\R^6_1$ --- $Sl(2,\H)$ is a double cover of the identity component\Footnote{%
Using a basis of quaternionic hermitian forms, it is an unpleasant but
straightforward calculation to establish a Lie algebra isomorphism
${\fk sl}(2,\H)\leftrightarrow{\fk o}_1(6)$.} of $SO_1(6)$, which itself
is isomorphic to the group of orientation preserving M\"obius transformations
of $S^4$.
Thus, restricting our attention to Euclidean 4-space $\{e_1h+e_2\,|\,h\in\H\}$,
the orientation preserving M\"obius transformations appear as fractional linear
transformations (cf.\Study,\Wilker)
$$\left(\matrix{h\cr1\cr}\right)\mapsto\left(\matrix{a_{11}&a_{12}\cr
        a_{21}&a_{22}\cr}\right)\left(\matrix{h\cr1\cr}\right)\simeq
        \left(\matrix{(a_{11}h+a_{12})(a_{21}h+a_{22})^{-1}\cr1\cr}\right).$$
If $s\neq0$ lies in the light cone of $\R^6_1$, $\langle s,s\rangle=0$, then
the corresponding quadratic form $v\mapsto s(v,v)$ annihilates exactly {\it one} direction
$v\H\subset\H^2$: $0=s(v,v)$ vanishes iff $0=|s_{11}v_1+s_{12}v_2|^2$ or
$0=|s_{21}v_1+s_{22}v_2|^2$ since at least one, $s_{11}$ or $s_{22}$ does
not vanish. Hence, we can identify a point $p=v\H\in\HP^1$ of the quaternionic
projective line --- the 4-sphere --- with the null line of quaternionic
hermitian forms in the Minkowski $\R^6_1$ that annihilate this point.
In homogeneous coordinates, this identification can be given by\Footnote{%
Note the analogy with the Veronese embedding.}
$$v=\left(\matrix{v_1\cr v_2\cr}\right)\leftrightarrow
  \left(\matrix{|v_2|^2&-v_1\bar{v}_2\cr-v_2\bar{v}_1&|v_1|^2\cr}\right)=s_v.
  \Eqno\Label\LightCone$$
Note, that with this identification, $\langle s_v,s\rangle=-s(v,v)$ for any
quaternionic hermitian form $s\in\R^6_1$. If $s=s_w$ is an isotropic form too,
then $\langle s_v,s_w\rangle=-\SD(v,w)$.

If, on the other hand, $\langle s,s\rangle=1$ we obtain --- depending on
whether $s_{11}=0$ or $s_{11}\neq0$ in the chosen basis $(e_1,e_2)$ of $\H^2$
--- $$s=\left(\matrix{0&-n\cr-\bar{n}&2d\cr}\right)\hskip2em{\rm or}\hskip2em
      s={1\over r}\left(\matrix{1&-m\cr-\bar{m}&|m|^2-r^2\cr}\right)$$
with suitable $n$ resp.\ $m\in\H$ and $d$ resp.\ $r\in\R$: the null cone of $s$
is a plane with unit normal $n$ and distance $d$ from the origin or a sphere
with center $m$ and radius $r$ in Euclidean 4-space $\{e_1h+e_2\,|\,h\in\H\}$.
Consequently, we identify the Lorentz sphere $S^5_1\subset\R^6_1$ with the
space of spheres and planes in Euclidean 4-space, or with the space of spheres
in $S^4$ --- as the readers familiar with the classical model (cf.\Blaschke) of
M\"obius geometry might already have suspected. The incidence of a point $p\in
S^4\cong\HP^1$ and a sphere $s\subset S^4$, i.e.\ $s\in S^5_1$, is equivalent
to $s(p,p)=0$ in our quaternionic model.
A key concept in

\Section{M\"obius differential geometry}
is that of (hyper-) sphere congruences and envelopes of sphere congruences:

\Definition{}{An immersion $f:M\to S^4$ is called an envelope of a
hypersphere congruence $s:M\to S^5_1$ if, at each point $p\in M$,
$f$ touches the corresponding sphere $s(p)$:
$f(p)\in s(p)$ and $d_pf(T_pM)\subset T_{f(p)}s(p)$.}

\noindent
According to our previous discussion, the first condition --- the incidence of
$f(p)$ and the corresponding sphere $s(p)$ --- is equivalent to $s(f,f)=0$ in
our quaternionic model. Calculating, for a moment, in a Euclidean setting
--- i.e.\ $s={1\over r}\left(\matrix{1&-m\cr-\bar{m}&|m|^2-r^2}\right)$ ---
we find $s(f,df)+s(df,f)={2\over r}(f-m)\cdot df$. Thus\Footnote{%
Note, that with the identification \LightCone\ of points in $\HP^1$ with
isotropic quaternionic hermitian forms, $s(f,df)+s(df,f)=-\langle s,df\rangle$
which gives the link with the classical model of M\"obius geometry.},

\Lemma{}{An immersion $f:M\to\HP^1$ is an envelope of a sphere congruence
$s:M\to S^5_1$ if and only if $s(f,f)=0$ and $s(f,df)+s(df,f)=0$.}

\noindent
Before going on, we introduce the symmetric space of point pairs:
given two (distinct) points of the quaternionic projective line $\HP^1$,
we may identify these points with a quaternionic linear transformation $P$
which maps a (fixed) basis $(e_1,e_2)$ of $\H^2$ to their homogeneous
coordinates --- or, in coordinates, with a matrix having for columns the
homogeneous coordinates of the two points. This linear transformation $P$
is obviously not uniquely determined by the two points in $\HP^1$: any
gauge transform $P\cdot H$ of $P$ with $H$ in the isotropy subgroup
$K:=\{H\in Gl(2,\H)\,|\,He_1=e_1\lambda,He_2=e_2\mu\}$ determines the
same point pair. Thus, the space $\P$ of point pairs in the conformal
4-sphere $\HP^1$ is a homogeneous space, $\P=Gl(2,\H)/K$.
Moreover, the decomposition ${\fk gl}(2,\H)={\fk k}\oplus{\fk p}$
with $$\matrix{
        {\fk k}&=&\{X\in{\fk gl}(2,\H)\,|\,Xe_1=e_1\lambda,Xe_2=e_2\mu\}\cr
        {\fk p}&=&\{X\in{\fk gl}(2,\H)\,|\,Xe_1=e_2\lambda,Xe_2=e_1\mu\}\cr}
        \Eqno\Label\Cartan$$
is a Cartan decomposition since $[{\fk k},{\fk k}]\subset{\fk k}$
$[{\fk k},{\fk p}]\subset{\fk p}$ and $[{\fk p},{\fk p}]\subset{\fk k}$
so that $\P$ is, in fact, a symmetric space.

Now, if $F=(f,\hat{f}):M\to Gl(2,\H)$ is a framing (lift) of a point
pair map $M\to\P$, a simple calculation using \LightCone\ shows that
$$Ff=\left(\matrix{0&0\cr0&1\cr}\right)\hskip1em\hbox{and}\hskip1em
        F\hat{f}=\left(\matrix{1&0\cr0&0\cr}\right)$$
if the relative scaling of $f$ and $\hat{f}$ is chosen such that $F$
takes values in the special linear group $Sl(2,\H)$.
Since $Sl(2,\H)$ acts by isometries on the space $\R^6_1$ of quaternionic
hermitian forms, for any sphere congruence $s:M\to S^5_1$ containing the
points of $f$ and $\hat{f}$, we have
$$Fs=\left(\matrix{0&s_0\cr\bar{s}_0&0\cr}\right)$$
with a suitable function $s_0:M\to S^3\subset\H$ taking values in the
unit quaternions. Passing to another set of homogeneous coordinates
by means of a gauge transformation
$(f,\hat{f})\mapsto(f\lambda,\hat{f}\hat{\lambda})$ results in
$s_0\mapsto\bar{\lambda}s_0\hat{\lambda}$. Thus, depending on a given
sphere congruence $s$, we may fix the homogeneous coordinates of $f$ and
$\hat{f}$ such that $s_0\equiv1$ --- leaving us with a scaling freedom
$(f,\hat{f})\mapsto(f\lambda,\hat{f}\bar{\lambda}^{-1})$ with $\lambda\in\H$.
A second sphere congruence $\tilde{s}$ (orthogonal to the first one) can be
used to further fix the scalings via $\tilde{s}_0\equiv i$ up to $\lambda\in\C$.
Giving a complete set of four accompaning orthogonal sphere congruences (by
fixing a third one, $\hat{s}$, to satisfy $\hat{s}_0\equiv j$) leaves us
with the familiar real scaling freedom, $\lambda\in\R$ (cf.\Blaschke).

Writing down the derivatives $df=f\varphi+\hat{f}\psi$ and $d\hat{f}=
f\hat{\psi}+\hat{f}\hat{\varphi}$ of $f$ and $\hat{f}$, we obtain the
connection form
$$\Phi:=F^{-1}dF
        =\left(\matrix{\varphi&\hat{\psi}\cr\psi&\hat{\varphi}\cr}\right)
        :TM\to{\fk gl}(2,\H)$$
of a framing $F:M\to Gl(2,\H)$.
A gauge transformation $(f,\hat{f})\mapsto(f\lambda,\hat{f}\hat{\lambda})$
of the frame will result in a change
$$\left(\matrix{\varphi&\hat{\psi}\cr\psi&\hat{\varphi}\cr}\right)\mapsto
        \left(\matrix{\lambda^{-1}\varphi\lambda&
        \lambda^{-1}\hat{\psi}\hat{\lambda}\cr
        \hat{\lambda}^{-1}\psi\lambda&
        \hat{\lambda}^{-1}\hat{\varphi}\hat{\lambda}\cr}\right)+
        \left(\matrix{\lambda^{-1}d\lambda&0\cr0&
        \hat{\lambda}^{-1}d\hat{\lambda}\cr}\right)\Eqno\Label\Gauge$$
of the connection form $\Phi$.
The integrability conditions $0=d^2\!f=d^2\!\hat{f}$ yield the Maurer-Cartan
equation $0=d\Phi+\Phi\wedge\Phi$ for the connection form:
the Gauss-Ricci equations for $f$ resp.\ $\hat{f}$, $$\matrix{
  0&=&d\varphi+\varphi\wedge\varphi+\hat{\psi}\wedge\psi\hfill\cr
  0&=&d\hat{\varphi}+\hat{\varphi}\wedge\hat{\varphi}+\psi\wedge\hat{\psi},\cr
  }\Eqno\Label\Gauss$$
and the Codazzi equations, $$\matrix{
  0&=&d\psi+\psi\wedge\varphi+\hat{\varphi}\wedge\psi\hfill\cr
  0&=&d\hat{\psi}+\hat{\psi}\wedge\hat{\varphi}+\varphi\wedge\hat{\psi}.\cr
  }\Eqno\Label\Codazzi$$
Note, that since the quaternions are not commutative, generally
$\varphi\wedge\varphi\neq0$. Moreover,
$d(\lambda\varphi)=d\lambda\wedge\varphi+\lambda d\varphi$,
$d(\varphi\lambda)=d\varphi\,\lambda-\varphi\wedge d\lambda$ and
$\overline{\varphi\wedge\psi}=-\bar{\psi}\wedge\bar{\varphi}$
for any quaternion valued 1-forms $\varphi$ and $\psi$ and function
$\lambda:M\to\H$.

If $s:M\to\R^6_1$ is a map into the vector space of quaternionic hermitian
forms, then its derivative, $ds:TM\to\R^6_1$ is a 1-form with values in
the quaternionic hermitian forms. If $Fs\equiv const$, this derivative can
be expressed in terms of the connection form $\Phi$ of $F$:
since $d(Fs)=0$, $$F\,ds=-F[s(.,\Phi)+s(\Phi,.)]\simeq
        -[Fs\cdot\Phi+\Phi^{\ast}\cdot Fs]\Eqno\Label\Derivative$$
when using the matrix representation for quaternionic hermitian forms.

\Section{Curved flats and Isothermic surfaces}
The concept of curved flats in symmetric spaces was first introduced by
D.\ Ferus and F.\ Pedit \FePe. In \BJPP\ it was then applied to the geometry
of isothermic surfaces in 3-space. To demonstrate our quaternionic calculus
at work, we are going to discuss curved flats in the symmetric space $\P$ of
point pairs in $\HP^1$. These will turn out to be Darboux pairs of isothermic
surfaces in 4-space: given a point pair map $(f,\hat{f}):M\to\P$, we choose
a framing $F:M\to Sl(2,\H)$ and write its connection form $\Phi=\Phi_{\fk
k}+\Phi_{\fk p}:TM\to{\fk sl}(2,\H)={\fk k}\oplus{\fk p}$. Then,

\Definition{}{A map $(f,\hat{f}):M\to\P$ into the symmetric space of point
pairs is called a curved flat if $\Phi_{\fk p}\wedge\Phi_{\fk p}=0$.}

\noindent
Note, that the defining equation is invariant under gauge transformations
\Gauge\ of $F$, i.e.\ does not depend on a choice of homogeneous coordinates.
Thus, the notion of a curved flat is a well defined notion
for a point pair map $(f,\hat{f}):M\to\P$.

In order to understand the geometry of a curved flat $(f,\hat{f}):M^2\to\P$
in the symmetric space of point pairs we will first express its connection form
in a simpler form, and then interpret it geometrically in a second step. We
start with an $Sl(2,\H)$-framing $F:M^2\to Sl(2,\H)$ and write its connection
form
$$\Phi=\left(\matrix{
        \varphi_1+\varphi_2j&\hat{\psi}_1+\hat{\psi}\,j\cr
        \psi_1+\psi\,j&\hat{\varphi}_1+\hat{\varphi}_2j\cr}\right)$$
in terms of complex valued 1-forms.
Using a rescaling $(f,\hat{f})\mapsto(f\lambda,\hat{f}\hat{\lambda})$ we
can achieve $\psi_1=0$; then, the curved flat equations read (we assume
$\psi\neq0$) $\hat{\psi}_1=0$ and $\hat{\psi}\wedge\bar{\psi}=0$.
A rescaling $(f,\hat{f})\mapsto(f\bar{\lambda},\hat{f}\lambda^{-1})$ with a
complex valued function $\lambda$ results in $(\psi,\hat{\psi})\mapsto
(\lambda^2\psi,\bar{\lambda}^{-2}\hat{\psi})$; as any 1-form on $M^2$ has
an integrating factor, we may assume $d\psi=0$, i.e.\ $\psi=dw$. Since
$\hat{\psi}\wedge\bar{\psi}=0$, $\hat{\psi}=\bar{a}^4d\bar{w}$ with a suitable
function $a:M\to\C$. From the Codazzi equations, $da\wedge dw=0$ --- thus,
by a holomorphic change $z_w=a^2$ of coordinates, $\psi=a^{-2}dz$ and
$\hat{\psi}=\bar{a}^2d\bar{z}$, or, after rescaling again with $\lambda=a$,
$\psi=dz$ and $\hat{\psi}=d\bar{z}$. Now, the Codazzi equations also yield
$\hat{\varphi}_2\wedge dz=\bar{\varphi}_2\wedge d\bar{z}$ and
$\hat{\varphi}_2\wedge d\bar{z}=\bar{\varphi}_2\wedge dz$.
Thus, $\varphi_2=q_1dz-\bar{q}_2d\bar{z}$ and
$\hat{\varphi}_2=-\bar{q}_1dz+q_2d\bar{z}$
with suitable functions $q_1,q_2:M\to\C$. This way,
$\varphi_2\wedge\bar{\varphi}_2=\hat{\varphi}_2\wedge\bar{\hat{\varphi}}_2$
such that $d\varphi_1=d\hat{\varphi}_1$ from the Gauss-Ricci equations.
With the ansatz $\hat{\varphi}_1-\varphi_1=2a$, we find that a rescaling
$(f,\hat{f})\mapsto(f\lambda,\hat{f}\lambda^{-1})$ with $\lambda=e^a$
yields $\varphi_1=\hat{\varphi}_1$.
At the same time, $(\psi,\hat{\psi})\mapsto(e^u\psi,e^{-u}\hat{\psi})$ with
$u=a+\bar{a}$. So, we end up with a connection form $$\Phi=\left(\matrix{
        i\eta+(q_1dz-\bar{q}_2d\bar{z})j&e^{-u}d\bar{z}\,j\cr
        e^udz\,j&i\eta+(-\bar{q}_1dz+q_2d\bar{z})j\cr}\right)
        \Eqno\Label\CFConnect$$
where $u:M\to\R$, $q_1,q_2:M\to\C$ and $\eta:TM\to\R$ is a real valued 1-form
--- remember that we have chosen an $Sl(2,\H)$-framing from the beginning.

In order to interpret this connection form geometrically, we first note that
all sphere congruences
$$s_c:=F^{-1}\left(\matrix{0&c\cr\bar{c}&0\cr}\right):M\to S^5_1$$
with $c=e^{i\vartheta}$ are enveloped by the two maps $f$ and $\hat{f}$:
$$Fds_c=-\left(\matrix{0&2[-\Re(\bar{c}q_1)dz+\Re(cq_2)d\bar{z}]j\cr
        2[\Re(\bar{c}q_1)dz-\Re(cq_2)d\bar{z}]j&0\cr}\right)$$
Thus, in the $\R^6_1$-model of M\"obius geometry, the $s_c$ can be viewed
as common normal fields of $f$ and $\hat{f}$. Using the identification
\LightCone\ of points in $\HP^1$ and isotropic lines in $\R^6_1$, we obtain
$$df=F^{-1}\left(\matrix{0&e^udz\,j\cr-e^udz\,j&0\cr}\right)
        \hskip1em\hbox{and}\hskip1em
  d\hat{f}=F^{-1}\left(\matrix{0&-e^{-u}d\bar{z}\,j\cr
        e^{-u}d\bar{z}\,j&0\cr}\right)$$
as the derivatives \Derivative\ of $f$ and $\hat{f}$. Calculating the
induced metrics
$$\langle df,df\rangle=e^{2u}|dz|^2\hskip1em\hbox{and}\hskip1em
        \langle d\hat{f},d\hat{f}\rangle=e^{-2u}|dz|^2$$
of $f$ and $\hat{f}$, and their second fundamental forms with respect to $s_c$,
$$\matrix{-\langle df,ds_c\rangle&=&\hfill
        e^u[-2\Re(\bar{c}q_1)|dz|^2+\Re(cq_2)(dz^2+d\bar{z}^2)],\cr
        -\langle d\hat{f},ds_c\rangle&=&
        e^{-u}[-2\Re(cq_2)|dz|^2+\Re(\bar{c}q_1)(dz^2+d\bar{z}^2)],\cr}$$
we see that $f$ and $\hat{f}$ have well defined principal curvature
directions (independent of the normal direction $s_c$) which do correspond
on both surfaces ($\{s_c\,|\,c\in S^1\}$ is a ``Ribaucour sphere pencil''),
and that $f$ and $\hat{f}$ induce conformally equivalent metrics on $M$.
Moreover, $z:M\to\C$ are conformal curvature line coordinates on both surfaces,
i.e.\ both surfaces are isothermic. Consequently, $(f,\hat{f}):M\to\P$
is a ``Darboux pair'' of isothermic surfaces in 4-space\Footnote{This
geometric description of Darboux pairs of isothermic surfaces can obviously
be used to define isothermic surfaces and Darboux pairs of any codimension
--- as the one below for Christoffel pairs can (cf.\Palmer). Note, that the
flatness of the normal bundle of a surface --- which is necessary to make
sense of the notion of curvature lines --- is a conformal notion, i.e.\ it
is invariant under conformal changes of the ambient space's metric.}:

\Definition{}{Two surfaces are said to form a Darboux pair if they envelope
a (nontrivial) congruence of 2-spheres (two orthogonal congruences of 3-spheres
in 4-space) such that the curvature lines on both surfaces correspond and the
induced metrics in corresponding points are conformally equivalent.}

\noindent
Conversely, if $(f,\hat{f}):M\to\P$ envelope two congruences of orthogonal
spheres, say $s_1,s_i:M\to S^5_1$, then the connection form
$$\Phi=\left(\matrix{\varphi_1+\varphi_2j&\hat{\psi}\,j\cr
                     \psi\,j&\hat{\varphi}_1+\hat{\varphi}_2j\cr}\right)$$
with complex 1-forms $\psi,\hat{\psi}:TM\to\C$. Assuming the curvature lines of
$f$ and $\hat{f}$ to correspond, and their induced metrics to be conformally
equivalent, we can introduce common curvature line coordinates:
$\psi=e^u\omega$ and $\hat{\psi}=e^{-u}\omega$, or
$\hat{\psi}=e^{-u}\bar{\omega}$.
In both cases, from the Gauss-Ricci equations
$\Re[d(\varphi_1-\hat{\varphi}_1)]=0$, so that after a suitable real
rescaling of $f$ and $\hat{f}$, $\Re(\varphi_1-\hat{\varphi}_1)=0$.
Then, in the first case, the Codazzi equations imply $u\equiv const$:
the sphere congruences enveloped by $f$ and $\hat{f}$ lie in a fixed linear
complex, consequently $f$ and $\hat{f}$ are congruent in some space of constant
curvature (cf.\Blaschke, \BJPP) --- and are not considered to form a Darboux
pair.
In the other case, the Codazzi equations yield $d\omega=0$ --- we have {\it
conformal} curvature line parameters, i.e.\ $f$ and $\hat{f}$ are isothermic;
we could also have concluded this from the fact that $f$ and $\hat{f}$
obviously form a curved flat:

\Theorem{}{A surface pair $(f,\hat{f}):M^2\to\P$ is a curved flat if and only
if $f$ and $\hat{f}$ form a Darboux pair.\hfil\break
In particular, two surfaces forming a Darboux pair are isothermic.}

\noindent
The $\fk k$-part --- see \Cartan\ --- of the Maurer-Cartan equation of a
$Gl(2,\H)$-framing reads $0=d\Phi_{\fk k}+\Phi_{\fk k}\wedge\Phi_{\fk k}+
\Phi_{\fk p}\wedge\Phi_{\fk p}$.
Thus, for a curved flat, $\Phi_{\fk k}=H^{-1}dH$ with a suitable $H:M\to K$:
if $\lambda$ and $\hat{\lambda}$ are given by
$$\lambda^{-1}d\lambda=i\eta+(q_1dz-\bar{q}_2d\bar{z})j\hskip1em\hbox{and}
  \hskip1em\hat{\lambda}^{-1}d\hat{\lambda}=i\eta+(-\bar{q}_1dz+q_2d\bar{z})j$$
then a gauge transformation
$(f,\hat{f})\mapsto(f\lambda^{-1},\hat{f}\hat{\lambda}^{-1})$
of our previous framing with connection form \CFConnect\ leaves us with
$$\Phi=\left(\matrix{0&\lambda(e^{-u}d\bar{z}\,j)\hat{\lambda}^{-1}\cr
        \hat{\lambda}(e^udz\,j)\lambda^{-1}&0\cr}\right)=:
        \left(\matrix{0&\hat{\omega}\cr\omega&0\cr}\right).$$
The Codazzi equations for this new framing simply read $d\omega=d\hat{\omega}
=0$ showing that $\bar{\omega}=df_0$ and $\hat{\omega}=d\hat{f}_0$ with
suitable maps $f_0,\hat{f}_0:M\to\H$. Here, we identify the two copies of
the quaternions sitting in ${\fk p}=\H\oplus\H$ as the eigenspaces of
${\rm ad}_C:{\fk p}\to{\fk p}$, $C={\left(\matrix{1&\hfill0\cr0&-1\cr}\right)}$,
by means of the real endomorphism $X\mapsto X^{\ast}$ of ${\fk p}$.
Note, that since the 1-forms
$\lambda^{-1}d\lambda,\hat{\lambda}^{-1}d\hat{\lambda}:TM\to\Im\H$
take values in the imaginary quaternions, $|\lambda|=|\hat{\lambda}|\equiv1$.
Consequently, the induced metrics of $f_0:M\to\H$ and $\hat{f}_0:M\to\H$,
$\H\cong\R^4$ considered as a Euclidean space, are
$$df_0\cdot df_0=e^{2u}|dz|^2\hskip1em\hbox{and}\hskip1em
        d\hat{f}_0\cdot d\hat{f}_0=e^{-2u}|dz|^2.$$
Moreover, with the common unit normal fields $n_c=-\lambda c\hat{\lambda}^{-1}$
of $f_0$ and $\hat{f}_0$, where $c=e^{i\vartheta}$, their second fundamental
forms become
$$\matrix{-df_0\cdot dn_c&=&\hfill
        e^u[-2\Re(cq_1)|dz|^2+\Re(\bar{c}q_2)(dz^2+d\bar{z}^2)],\cr
        -d\hat{f}_0\cdot d\hat{n}_c&=&
        e^{-u}[-2\Re(\bar{c}q_2)|dz|^2+\Re(cq_1)(dz^2+d\bar{z}^2)].\cr}
        \Eqno\Label\ChristoffelII$$
Thus, $f_0$ and $\hat{f}_0$ are two isothermic surfaces that carry common
curvature line coordinates --- and, $\hat{f}_0$ and $\bar{f}_0$ have parallel
tangent planes.
Hence, we define\Footnote{If $f_0,\hat{f}_0:M^2\to\Im\H$, this definition yields
the classical notion of a Christoffel pair (cf.\BJPP).}:

\Definition{}{Two (non homothetic) surfaces $f_0,\hat{f}_0:M^2\to\H$ with
parallel tangent planes in corresponding points are said to form a Christoffel
pair if the curvature lines on both surfaces correspond and the induced metrics
are conformally equivalent.}

\noindent
Conversely, if two surfaces $f_0,\hat{f}_0:M^2\to\H$ carry conformally
equivalent metrics and have parallel tangent planes in corresponding points
$f_0(p)$ and $\hat{f}_0(p)$ then\Footnote{If $p$ is not an umbilic for either
surface, it follows that the principal curvature directions of both surfaces
correspond. In case one of the surfaces is totally umbilic we need also to
assume that the curvature lines on both surfaces coincide --- otherwise we
might find two associated minimal surfaces.},
$df_0=\lambda e^u\psi\,j\hat{\lambda}^{-1}$ and
$d\hat{f}_0=\pm\lambda e^{-u}\psi\,j\hat{\lambda}^{-1}$, or
$d\hat{f}_0=\lambda e^{-u}\bar{\psi}\,j\hat{\lambda}^{-1}$
with a real valued function $u$, a complex 1-form $\psi:TM\to\C$ and
suitable quaternionic functions $\lambda,\hat{\lambda}:M\to\H$ --- where
$|\lambda|=|\hat{\lambda}|\equiv1$ without loss of generality.
In the first case, the integrability conditions yield $0=du\wedge\psi$
showing that $u\equiv const$. Consequently, $\hat{f}_0$ is homothetic to $f_0$
--- and $f_0$ and $\hat{f}_0$ are not considered to form a Christoffel pair.
In the second case, $d\bar{f}_0\wedge d\hat{f}_0=d\hat{f}_0\wedge d\bar{f}_0=0$.
Hence, the surface pair $f_0,\hat{f}_0:M\to\H$ gives rise to a curved flat
by integrating $\Phi:=\left(\matrix{0&d\hat{f}_0\cr d\bar{f}_0&0\cr}\right)$
--- we obtain the following

\Theorem{}{Two surfaces $f_0,\hat{f}_0:M^2\to\H$ form a Christoffel pair if
and only if
$d\bar{f}_0\wedge d\hat{f}_0=d\hat{f}_0\wedge d\bar{f}_0=0$.\hfil\break
In particular, two surfaces forming a Christoffel pair are isothermic.}

\noindent
Curved flats --- or, Darboux pairs of isothermic surfaces --- naturally arise
in 1-parameter families \FePe: if $\Phi=\Phi_{\fk k}+\Phi_{\fk p}$ denotes one
of the connection forms associated to a curved flat $(f,\hat{f}):M^2\to\P$,
then, with a real parameter $\varrho\in\R$, all the connection forms
$$\Phi_{\varrho}:=\Phi_{\fk k}+\varrho^2\Phi_{\fk p}:
  TM^2\to{\fk sl}(2,\H)={\fk k}\oplus{\fk p}\Eqno\Label\Family$$
are integrable and give rise to curved flats
$(f_{\varrho},\hat{f}_{\varrho}):M^2\to\P$;
in fact, if the connection forms \Family\ are integrable for more than one
value of $\varrho^2$, then the associated point pair maps are necessarily
curved flats.
From \Gauge, we learn that this 1-parameter family of curved flats does not
depend on the framing chosen to describe the curved flat $(f,\hat{f})$.
Moreover, sending the parameter $\varrho\to0$, and rescaling
$(f_{\varrho},\hat{f}_{\varrho})\mapsto
  (\varrho^{-1}f_{\varrho},\varrho\hat{f}_{\varrho})$
or
$(f_{\varrho},\hat{f}_{\varrho})\mapsto
  (\varrho f_{\varrho},\varrho^{-1}\hat{f}_{\varrho})$
at the same time, provides us with
$$(f_{\varrho=0},\hat{f}_{\varrho=0})
=\left(\matrix{1&0\cr\bar{f}_0&1\cr}\right)
\hskip1em\hbox{or}\hskip1em
(f_{\varrho=0},\hat{f}_{\varrho=0})
=\left(\matrix{1&\hat{f}_0\cr0&1\cr}\right).$$
Hence, we may think of the Christoffel pair $(f_0,\hat{f}_0)$ --- that is,
as before, associated to a 1-parameter family of curved flats by integrating
$$\Phi_{\varrho}
  =\left(\matrix{0&\varrho^2d\hat{f}_0\cr\varrho^2d\bar{f}_0&0\cr}\right)$$
--- as a limiting case for the Darboux pairs $(f_{\varrho},\hat{f}_{\varrho})$.
Comparison with \Gauge\ shows that the spectral parameter $\varrho$ corresponds
to the scaling ambiguity of the members of a Christoffel pair: one of the
surfaces of a Christoffel pair is determined by the other only up to a
homothety (and translation).

We will use those facts to discuss perturbation methods (cf.\UmYa) for the
construction of constant mean curvature surfaces and, in particular, for
Bryant's Weierstrass type representation \Bryant\ for

\Section{Constant mean curvature surfaces}
in hyperbolic space forms. We restrict our attention to codimension 1
by assuming that our surfaces lie in a fixed conformal 3-sphere, say $s_1$.
Thus the connection form \CFConnect\ of a Darboux pair $(f,\hat{f}):
M^2\to\HP^1$ takes the form
$$\Phi=\left(\matrix{
  i[\eta+{1\over2}(e^uHdz-e^{-u}\hat{H}d\bar{z})j]&e^{-u}d\bar{z}\,j\cr
  e^udz\,j&i[\eta+{1\over2}(e^uHdz-e^{-u}\hat{H}d\bar{z})j]\cr}\right)
  \Eqno\Label\CodimenOne$$
where the (real) functions $H,\hat{H}$ can be interpreted as the mean
curvature functions of the members $f_0$ and $\hat{f}_0$ of the limiting
Christoffel pair: from \CodimenOne\ we see that a rescaling $(f,\hat{f})
\mapsto(f\lambda,\hat{f}\lambda)$ will provide us with $\Phi_{\fk k}=0$,
such that $df_0,d\hat{f}_0:TM\to\Im\H$. The second fundamental forms
\ChristoffelII\ with respect to the remaining common normal field
$n_i=-\lambda i\lambda^{-1}$ become
$$\matrix{-df_0\cdot dn_i&=&\hfill
        He^{2u}|dz|^2-{1\over2}\hat{H}(dz^2+d\bar{z}^2)],\cr
        -d\hat{f}_0\cdot d\hat{n}_i&=&
        \hat{H}e^{-2u}|dz|^2-{1\over2}H(dz^2+d\bar{z}^2)].\cr}$$
The Codazzi equations \Codazzi\ yield $\eta={i\over2}(-u_zdz+u_{\bar{z}}
d\bar{z})$ and from \Gauss\ we recover the classical Gauss equation
$0=u_{z\bar{z}}+{1\over4}(H^2e^{2u}-\hat{H}^2e^{-2u})$
holding for both surfaces $f_0$ and $\hat{f}_0$, and the classical Codazzi
equations
$dH\wedge e^udz=d\hat{H}\wedge e^{-u}d\bar{z}.$
Hence, $H\equiv const$ if and only if $\hat{H}\equiv const$, reflecting the
fact that a pair of parallel constant mean curvature surfaces, or a minimal
surface and its Gauss map form Christoffel pairs (cf.\JePe).

Calculating the derivative of the sphere congruence $s_i$ enveloped by the
two surfaces $f$ and $\hat{f}$ --- which form the Darboux pair associated
with the Christoffel pair $(f_0,\hat{f}_0)$ --- we find
$$Fds_i=\left(\matrix{0&(He^udz-\hat{H}e^{-u}d\bar{z})j\cr
        (-He^udz+\hat{H}e^{-u}d\bar{z})j&0\cr}\right)
        =H\cdot F\,df+\hat{H}\cdot F\,d\hat{f}.$$
Hence, the vector $N:=s_i-Hf-\hat{H}\hat{f}$ is constant as soon as one of the
mean curvatures, $H$ or $\hat{H}$, is.
In order to interpret this fact geometrically, we have to distinguish two cases:

If $H\hat{H}\neq0$, i.e.\ $(f_0,\hat{f}_0)$ is equivalent to a pair of parallel
constant mean curvature surfaces, $\langle N,{2\over\hat{H}}f\rangle\equiv1$
and $\langle N,{2\over H}\hat{f}\rangle\equiv1$. Consequently (cf.\Blaschke),
the two surfaces ${1\over\hat{H}}f,{1\over H}\hat{f}:M^2\to s_1\simeq S^3
\subset\HP^1$ can be interpreted as surfaces in the space $M^3_N:=\{y\in\R^6_1
\,|\,\langle N,y\rangle=1,\langle s_1,y\rangle=0\}$ of constant curvature
$\kappa=-\langle N,N\rangle=-(1-H\hat{H})$. Their induced metrics are
$$\matrix{
  \langle d({2\over\hat{H}}f),d({2\over\hat{H}}f)\rangle
        ={4\over\hat{H}^2}e^{2u}|dz|^2
    &\hskip.5em\hbox{and}\hskip.5em&
  \langle d({2\over H}\hat{f}),d({2\over H}\hat{f})\rangle
        ={4\over H^2}e^{-2u}|dz|^2\cr}$$
while, with the unit normal fields $t=s_i-{2\over\hat{H}}f$ and
$\hat{t}=s_i-{2\over H}\hat{f}$ in that space $M^3_N$, their second
fundamental forms become
$$\matrix{-\langle d({2\over\hat{H}}f),dt\rangle&=&\hfill
    {4\over\hat{H}^2}e^{2u}(1-{1\over2}H\hat{H})\,|dz|^2+(dz^2+d\bar{z}^2)\cr
          -\langle d({2\over H}\hat{f}),d\hat{t}\rangle&=&
    {4\over H^2}e^{-2u}(1-{1\over2}H\hat{H})\,|dz|^2+(dz^2+d\bar{z}^2)\cr}$$
--- showing that both surfaces have the same constant mean curvature
$1-{1\over2}H\hat{H}$.
As a special case, $H=1$ and $\hat{H}=2$, this provides the well known
relation between constant mean curvature surfaces in Euclidean space $\R^3$
and minimal surfaces in the 3-sphere $S^3$.

If $H\hat{H}=0$, one of the two surfaces, $f_0$ or $\hat{f}_0$ is a minimal
surface, say $\hat{H}=0$, while the other is homothetic to its Gauss map,
say $n=H\,f_0$. Now, the surface ${2\over H}\hat{f}:M^2\to M^3_N$ lies
in hyperbolic space, $\kappa=-1$, while $f$ is the hyperbolic Gauss map
(cf.\Bryant) of ${2\over H}\hat{f}$ since $\langle N,f\rangle\equiv0$,
i.e.\ $f$ takes values in the infinity boundary $N\in S^5_1$ of $M^3_N$.
As before, the mean curvature of ${2\over H}\hat{f}:M^2\to M^3_N$ is easily
calculated to be constant $=1$.
This is how Bryant's Weierstrass type representation \Bryant\ for surfaces
of constant mean curvature 1 in hyperbolic 3-space $H^3$ can be obtained
in this context:
we write the differential $d\hat{f}_0={1\over2}(i+gj)\bar{\omega}j(i+gj)$
of a minimal immersion $\hat{f}_0:M^2\to\R^3$ (and its Christoffel transform,
its Gauss map $f_0=(i+gj)i(i+gj)^{-1}:M^2\to S^2$) in terms of a holomorphic
1-form $\omega:TM^2\to\C$ and the stereographic projection $g:M\to\C$ of its
meromorphic Gauss map.
Then, the constant mean curvature surface $\hat{f}:M^2\to H^3$ (and its
hyperbolic Gauss map $f:M^2\to N\simeq S^2$) are obtained by integrating
the connection form\Footnote{With the ansatz $F=\left(\matrix{\hfill
        2(\overline{x_{21}g+x_{22}})(i+gj)^{-1}& j(x_{21}i-x_{22}j)\cr
        2j(\overline{x_{11}g+x_{12}})(i+gj)^{-1}&-(x_{11}i-x_{12}j)\cr
  }\right)$, the common form of Bryant's representation is obtained as
$xx^{\ast}:M^2\to H^3\cong\{y\in Gl(2,\C)\,|\,y=y^{\ast}\}$ where the
scalar product on $H^3$ is given by $|y|^2=-\det(y)$.}
$$\Phi=\left(\matrix{0&{1\over2}(i+gj)\bar{\omega}j(i+gj)\cr
        -2(i+gj)^{-1}dg\,j(i+gj)^{-1}&0\cr}\right),
        \Eqno\Label\Minimal$$
to the framing $(f,\hat{f})\simeq F:M^2\to Gl(2,\H)$ where $dF=F\Phi$.
In fact, introducing the spectral parameter \Family, surfaces of constant mean
curvature $c$ in hyperbolic space forms of curvature $\kappa=-c^2$ arise by
``perturbation'' of minimal surfaces in Euclidean 3-space (cf.\UmYa).

A closer look at the connection form \Minimal\ suggests that the classical
Weierstrass representation for minimal surfaces in $\R^3$ can be interpreted
as a Goursat type transformation of the plane:
considering $gj,\int\bar{\omega}j:M^2\to\C\,j$ as a (highly degenerate)
Christoffel pair, the corresponding minimal surface (and its Gauss map, {\it
its} Christoffel transform) is obtained as a Christoffel transform of a M\"obius
transform, the stereographic projection $f_0={1\over1+|g|^2}[(1-|g|^2)i+2gj]$,
of $gj$ (``the'' Christoffel transform of $\int\bar{\omega}j$).
This Goursat type transformation can (obviously) be generalized to arbitrary
Christoffel pairs of isothermic surfaces: if $f_0,\hat{f}_0:M^2\to\H$ form a
Christoffel pair, then, for any (constant) $a\in\H$, the quaternionic 1-forms
$(a+\bar{f}_0)^{-1}d\bar{f}_0(a+\bar{f}_0)^{-1}$ and
$(a+\bar{f}_0)d\hat{f}_0(a+\bar{f}_0)$ are closed
--- and consequently give rise to a new Christoffel pair.

\bigskip
\noindent{\bf Acknowledgements:} I would like to thank the many people who
contributed to this paper by questions and discussions.
In particular, I would like to thank the members of the GANG, my host at
GANG, Franz Pedit, Ulrich Pinkall who was visiting UMass during the winter
1995/96, and Ian McIntosh.

\References

\bye